\documentclass[12pt,a4paper]{article}
\usepackage[T2A]{fontenc}
\usepackage[cp1251]{inputenc}
\usepackage{amsmath,amssymb,amsthm,amsfonts,amscd}

 \theoremstyle{plain}
\newtheorem{theorem}{Theorem}
\newtheorem*{main*}{Main Theorem}
\newtheorem*{theorem*}{Theorem}
\newtheorem{lemma}[theorem]{Lemma}

\theoremstyle{definition}
\newtheorem{definition}[theorem]{Definition}

\newtheorem*{remark*}{Remark}
\newtheorem*{example*}{Example}

\def\Z{{\Bbb Z}}
\def\R{{\Bbb R}}

\def\N{{\Bbb N}}

\def\A{{\bf A}}

\def\F{{\bf F}}
\def\G{{\bf G}}

\def\L{{\bf L}}

\def\dist{{\operatorname{dist}}}
\def\grad{{\bf{grad}}}
\def\div{{\operatorname{div}}}
\def\rot{{\operatorname{rot}}}
\def\R{{\Bbb R}}
\def\B{{\bf B}}

\def\L{{\bf L}}
\def\valpha{\vec{\alpha}}

\date{}

\begin{document}
\sloppy

\title{On a higher analog of the asymptotic ergodic Hopf invariant}
\date{IZMIRAN, Russia, Moscow region}
\author{Petr M. Akhmet'ev \footnote{The author supported by Russian Foundation of Basic Research Grant No. 11-01-
00822.}}

\maketitle

\begin{abstract}
A particular results toward a positive solution of a problem by V.I.Arnol'd about a higher analog of the ergodic asymptotic invariant of magnetic fields is
presented.
\end{abstract}

\[  \]
Key words: magnetic field,  the asymptotic Hopf invariant, the Gauss integral, helicity, a finite order invariant of classical links,
the Birkhoff Theorem.
\[  \]
Codes MSC: 35Qxx, 57Mxx, 76Fxx.

\section{Introduction}

V.I.Arnol'd  formulated the following problem
\cite{Arn}, Problem 1984-12:
"To transform asymptotic ergodic definition of the Hopf invariant of divergence-free vector fields
to a theory by S.P.Novikov, which generalized Whitehead product in homotopy group of spheres"'.     

In the paper an asymptotic higher invariant is constructed using a generalized Massey invariant (Massey invariants
are determined by Whitehead integrals). We shall consider the 3-dimensional case, this case is the most important for applications. 

A divergence-free vector field in $\R^3$ with a compact support will be called a magnetic field. 
Magnetic fields form magnetic tubes, which are linked and knotted solid tori in the space. The goal is to define an invariant, which is preserved
with respect to volume-preserved transformations of the domain with a magnetic field, and which detects a geometrical complicity of the
magnetic tubes configuration. The asymptotic ergodic Hopf invariant, was proposed by V.I.Arnol'd in 1974 (see
\cite{Arn-Kh}) is generalized.  A new invariant is not a function on linking numbers (asymptotic ergodic) of
magnetic lines.

For magnetic fields with closed magnetic lines an asymptotic invariants has a combinatorial meaning: this invariants
are expressed using finite-order invariants of various multi-component links, which are formed by magnetic lines of the field.
The statement 1  of the main theorem  
$\ref{mainTh}$ (Existence) is particularly proved for generic magnetic fields, but the statement 2 (Invariance) is proved using this assumption.

In section 1 we recall the Main Theorem from \cite{A3}, which is based on Birkhoff Theorem. 
This  approach uses configuration spaces, associated with magnetic lines. Similar approach was developed in   
\cite{TGKMSV}.

In section 2 we recall required results from
\cite{A1}, where   "`a Milnor-Massey type"' invariant $M$ for magnetic fields inside 3 disjoint magnetic tubes is constructed.
 The higher invariant  $M$ is well-defined without additional assumptions, even in the case the pairwise integral linked numbers of the magnetic tubes are arbitrary.  

In section 3 the invariant $M$ is applied for magnetic tubes with identity Poincar$\acute{\rm{e}}$ recurrence mapping is investigated.
A limit, when the thickness of the magnetic tube vanishes, while magnetic flow reminds fixed is considered. 
The formula for the invariant $M$ in this limit is written as an integral over the corresponding configuration space, associated with the central lines of the magnetic tubes. This result is interesting by itself as a particular example to a general problem, investigated in \cite{C-K-Sh}.

In section 4 generic magnetic fields, represented by a finite number of magnetic tubes, are investigated.
The invariant $M$  as asymptotic ergodic invariant is defined, using results of section 1. Assuming that all magnetic lines in the tube are closed, or a magnetic field is structured in separated magnetic tubes, the invariant $M$ is expressed 
from the corresponding combinatorial invariant, which is derived in sections 3 and 4. 
For frozen-in magnetic fields with closed magnetic lines $M$ is an invariant with respect to volume-preserved diffeomorphisms of the domain.

This paper was presented at the International Conference in the Euler Institute (St-Petersburg) 10-14 September (2007),
and at the Rokhlin Memorial 11-16 January (2007) in the Euler Institute (St-Petersburg), and at the International Conference "Nonlinear Equations and Complex Analysis"
(Bashkortostan, Bannoe Lake) 13 - 17 December (2009), and at the workshop "`Entanglement and Linking"' (Pisa), 
18 May 2011 - 19 May (2011), and at the conference "`Quantized Flux in Tightly Knotted and Linked Systems"' (Cambridge) 3-7 December (2012). 
The author is grateful to  Prof. V.N.Vassiliev, O.Karpenkov, R. Komendarczyk, E.A.
Kudryavtseva, K.M.Kuzanyan, S.A.Melikhov, L.P.Plakhta, Prof. V.B.Semikoz, Prof. D.D.Sokolov,  Prof. A.B. Sossinsky, Prof. M.Spera and to Prof. R.L.Ricca 
 for discussions. 
 The author is grateful to the  participants of the seminar on analytic differential equations (MIRAN, Moscow) for discussion.

\section{Functionals of magnetic fields}

Let
$\B$, $\div(\B)=0$, is a smooth magnetic field in  $\R^3$ with a compact support 
$U \subset \R^3$, which is a manifold with a boundary. We assume that the magnetic field $\B$ is tangent
to the boundary  
$\partial U$ and there are no vanished inside $U$. In this case we say that a support $U$ is called a magnetic tube. 
All the $C^{\infty}$--magnetic fields in magnetic tubes are formed a space, which is denoted by $\Omega$.

\begin{definition}\label{local}

The configuration space $K_{q,r}$ is defined as following. Assume that a collection of $r$ magnetic lines
$L_1, \dots , L_r$ of the magnetic field $\B$, which is parametrized of the segments  $[0,T]$, started at the prescribed points
$\{l_1, \dots, l_r\}$ of the domain $\Omega$ correspondingly.
The subcollection $\{l_1; x_1, \dots, x_{q}\}$ of the full collection consists of  $q$ points, which are on the first magnetic line  $L_1$ of the magnetic field $\B$, the subcollection  
$\{l_2;x_{q+1}, \dots, x_{2q}\}$  consists of  $q$ points, each point belongs to the second
magnetic line $L_2$ of $\B$, e.t.c., the last subcollection
$\{l_r;x_{q(r-1)+1}, \dots, x_{qr}\}$ of the full collection consists of  $q$ points, each point belongs to the $r$-th magnetic line
$L_r$ of $\B$. Each point  $x_{qj+i}$ is well-defined by the time-variable $t_{qj+i}$, $1 \le j \le r$, $1 \le i \le q-1$, $0\le t_{qj+i} \le T$,
which is the time of the evolution of the point
$l_{j}$  into the point $x_{qj+i}$ by the magnetic flow. 

Let a real-valued functional (this definition makes sense for a multi-valued functional)  
 $\bar{I}: \Omega \to \R$ be well-defined. Let us say the functional $\bar{I}$ is of a finite-order, if it is defined as the average $\bar{I} = I(l_1,l_2,\dots,l_r)$ over the all collections $\{l_1, \dots, l_r\}$ of
the asymptotic limits for $T \to +\infty$ (called Cesaro averages) of integrals $\int f dx_{1,1}\dots dx_{rq}$ of a function $f: K_{q,r} \to \R$ over all finite collections  
 $\{l_1;x_1, \dots, x_{qr}\} \in K_{q,r}$ with fixed $\{l_1,\dots,l_r\}$. 
\end{definition}

Assume $\B \in \Omega$.  
Denote by $K_{r,q;T} \subset K_{r,q}$ a compact subspace in the configuration space, for which each time-variable coordinate belongs to the segment
 $[0,T]$. Let us formulate the definition of limiting tensor.

\begin{definition}\label{tenzor}
Assume a function
\begin{eqnarray}\label{Fq}
A: K_{r,q} \to \R
\end{eqnarray} 
is integrable on each subspace $K_{r,q;T} \subset K_{r,q}$.
Let us say that an integrable non-negative function
\begin{eqnarray}\label{a}
a^{[q]}: (U^q)^r \to \R_+,
\end{eqnarray}
which, possibly, tends to $+\infty$, when a point in the origin tends to the thick diagonal $Diag \subset (U^q)^r $,
is called a  {\it{limiting tensor}} for $(\ref{Fq})$, if there exists  $T_0 \ge 0$, such that for an arbitrary $T > T_0$ the function
$a^{[q]} \circ F_q: K(r,q) \to (U^q)^r \to \R_+ $ (this function is integrable, because the function  $(\ref{a})$ is integrable) and
$\vert f \vert : K_{r,q} \to \R_+$  satisfies the following equation:
\begin{eqnarray}\label{inneq}
\int   \vert A \vert dK_{r,q;T} \le \int a^{[q]} \circ F_q \quad dK_{r,q;T}.
\end{eqnarray}
\end{definition}
\[  \]
In the diagram above the evaluation mapping is used: 
\begin{eqnarray}\label{diag1}
K_{r,q}& \stackrel{F_q}{\longrightarrow} & ((U)^q)^r.
\end{eqnarray}
The evaluation mapping is used to investigate what's happening if we omit all the coordinates $\R_{1,i}, \R_{2,j}$ of points of the configuration spaces. 

The following statement is proved in \cite{A3}.

\begin{theorem}$\label{th4}$
Let $I(l_1,x_{1,1}, \dots, x_{r,q})$ is defined as a polynomial of functions $f_j \circ F_q(j): K_{r,1}(l_j,x_{j,1}, \dots, x_{j,q}) \to \R$,
$f_j: U^q \to \R$, $F_q(j): K_{r,1} \to U^q$.

--1. Assume that  for an arbitrary $j$ the function
$f_j: U^r \to \R$, $1 \le j \le q$ is integrable.

--2. 
Assume that  there exists a limiting tensor
$(\ref{a})$ for $I$ in the sense of Definition $\ref{tenzor}$. 
\[  \]

Then the asymptotic mean value $\bar{I}(l_1,\dots,l_r)$ of $I(l_1,x_{1,1}, \dots, x_{r,q})$ with respect to the coordinates $\R_{1,i}, \R_{2,j}$ is well-defined except, possibly, a  subset of $\{l_1, \dots l_r\}$ in
$U^r$ of zero measure, the function $\bar{I}(l_1,\dots,l_r): U^r \to \R$ is integrable (belongs to $L^1$) and invariant with respect to the magnetic flow on $U^r$.

\end{theorem}

\begin{definition}{\bf{Return Condition}}\label{vosvr}
Let $\B$, $\div(\B)=0$ be a smooth magnetic field with the support coincided with a magnetic tube  $U \subset \R^3$.
Let $\{g^t: U \to
U\}$ be the magnetic flow, generated by the magnetic field $\B$. Let us say that the magnetic field
$\B$ satisfies the Return Condition, if for an arbitrary point  $x \in \R^3$ there exists a real
$t_0>0$, $t_0=t_0(x)$, such that $g^{t_0}(x)=x$.
\end{definition}

 Return Condition implies that each magnetic line of $\B$ is closed
(we do not assume that the Poincar\'e mapping of a transversal section to itself is the identity). 
Evidently, there exist various examples of magnetic fields inside the only magnetic tube, which satisfy Definition 
$\ref{vosvr}$, and which have a finite number of ramifications of finite-orders. 
In this examples knots and links, which are formed by regular trajectories inside the tube are characterized by various invariants of isotopy classes. 
\[  \]

\section{Preliminaries: the formula for the invariant $M$ for a magnetic field inside 3 disjoint magnetic tubes}

Let us recall the integral formula from  \cite{A1}. Assume that a magnetic field $\B$ is localized into 3 disjoint arbitrary magnetic tubes
$$\B=\B_1 \cup \B_2 \cup \B_3. $$
assume that each magnetic tube is presented by a solid torus $U_i \subset \R^3$, in this torus a prescribed coordinate system 
$U_i = D^2 \times S^1$ is given by an arbitrary (volume-preserved) coordinate diffeomorphism. In this section we assume that the magnetic field $\B_i$, which is concentrated inside the magnetic tube $U_i$, is arbitrary. The magnetic flow of $\B_i$ trough the
transversal section of $U_i$ is parallel to the oriented central line $(0\times S^1)$ of the tube $U_i$.

Denote by
$\A_1, \A_2, \A_3$ the vector-potential of $\B$ in the corresponding magnetic tube: 
$$ \rot \A_i = \B_i, $$
denote by $(1,2),(2,3),(3,1)$ the integral linking coefficients, which are given by the integrals:
$$ (i,j) = \int_{U_i} (\A_j,\B_i) dU_i = \int_{U_j}(\A_i,\B_j) dU_j, \quad i,j = 1,2,3, \quad i \ne j. $$

Define an auxiliary vector-field
 $\F'$ in $\R^3$ by the formula: 
\begin{eqnarray}\label{FF}
\F' = (1,3)(2,3)\A_1 \times \A_2 + (2,1)(3,1)\A_2 \times \A_3 + (3,2)(1,2)\A_3 \times \A_1. 
\end{eqnarray}
Evidently, we get:
$$\div(\F') = (1,3)(2,3)[(\B_1,\A_2)-(\A_1,\B_2)] + \dots =$$
$$ (1,3)(2,3)(\B_1,\A_2)-(3,2)(1,2)(\A_3,\B_1) + \dots,$$
where $\dots$ means the terms, which are given by the cyclic permutation of the indexes. 

The restriction of the potential 
$(1,3)(2,3)\A_2 - (3,2)(1,2)\A_3$ on the magnetic tube  $U_1$ coincides with the gradient of a function, which is well-defined up to an additive constant. Such a function is denoted by
 $\phi_1 : U_1 \to \R$,  two analogous functions are denoted by 
$\phi_2 : U_2 \to \R$, $\phi_3 : U_3 \to \R$. These functions are defined below and is called the scalar potentials. 

Define the scalar potentials explicitly.
Denote by $\varphi_{j,i} : U_i \to \R^1$, $i \ne j$ a multivalued function, which is uniquely defined by the integration
of the vector-potential
$\A_j$ over the magnetic tube  $L_i$ along paths, which start at the marked point 
in $U_i$. This multivalued function is of the period $(i,j)$ is given by the formula:

\begin{eqnarray}\label{6}
\varphi_{j,i}(x_i) = \int_{pt_i}^{x_i} \A_j(x_i)  dx_i + C_{j,i}, \quad x_i
 \in U_i,
\end{eqnarray}
where $pt_i \in \Sigma_i \subset U_i$ is the marked point in $U_i$. Below we will consider only the case, when the magnetic tube is thin, and the surface $\Sigma$  coincides with a disk, which determines the transversal section of the magnetic tube. Assume that the coordinate system in   $U_i$ agrees with $\Sigma_i$, this means that
$U_i = \Sigma_i \times S^1$,  $\Sigma_i \cong D^2$. The coordinate lines along  $S^1$ coincide with the gradient of the prescribed multivalued functions 
$\varphi_{j,i}$, $j \ne i$.

By the construction $\varphi_{j,i}(pt_i)=C_{i,j}$, where the prescribed constant is additively well-defined up to $(i,j)k$, $k \in \Z$. 
The multivalued function $\varphi_{j,i}$ is in the formula  (19)\cite{A1}, this function is well-defined by the equation  (13)\cite{A1}.

Let us define the scalar potential $\phi_i: U_i \to \R^1$ by the formula: 
\begin{eqnarray}\label{8.1}
\phi_1 = (3,1) \varphi_{2,1} -  (1,2) \varphi_{3,1},
\end{eqnarray}
\begin{eqnarray}\label{8.2}
\phi_2 = (1,2) \varphi_{3,2} -  (2,3) \varphi_{1,2},
\end{eqnarray}
\begin{eqnarray}\label{8.3}
\phi_3 = (2,3) \varphi_{1,3} -  (3,3) \varphi_{3,1}.
\end{eqnarray}
The formulas $(\ref{8.1} )$--$(\ref{8.3})$ coincide with the formula  (16)-(18)\cite{A1}. 
The scalar potential $\phi_i$, $i=1,2,3$, is given by a scalar function, which is well defined up to an additive constant, let us call this constant 
a normalized constant. The normalized constant is defined below by the condition   $(\ref{mean})$.

The formulas 
$(\ref{8.1})$--$(\ref{8.3})$ should be investigated more carefully. The multivalued function 
$\varphi_{j,i}$ is decomposed into the sum of 3 terms. 
The first term, which is denoted by  $\varphi_{j,i}^0(x_i)$, $x_i \in U_i$, is a multivalued function, which is lifted into a linear function 
with respect to $\R$-coordinate on the universal covering   $\tilde U_i$ into the prescribed coordinate system  $\R \times \tilde \Sigma_i \cong \tilde U_i$,
and this function is zero on the surface $\tilde \Sigma_i$, the transversal section of the covering over the magnetic tube $\tilde U_i$. 

The second term is a function, which is denoted by
$\varphi_{j,i}^{var}: U_i \to \R$. The following condition is satisfied: 
\begin{eqnarray}\label{mean}
\int (\A^{0}_{j,i},\B_i)\varphi_{j,i}^{var} dU_i = 0,
\end{eqnarray}
where $\A^0_{j,i}  = \grad \varphi_{j,i}^0$. 
 
The third term is the constant, which is denoted by  
$C_{i,j}$, this constant determines the value of $\phi_i$ at the point  $pt_i$.
With respect to the definition above, we get:
\begin{eqnarray}\label{7}
\varphi_{j,i}(x_i) = \varphi_{j,i}^0 + \varphi_{j,i}^{var} +
C_{j,i}.
\end{eqnarray}
The terms in this formula depended on a coordinate system inside
$U_i$, but the integral $(\ref{MM})$ is well-defined below. 

Let us rewrite the formula ($\ref{8.1}$), ($\ref{8.2}$), ($\ref{8.3}$), using the normalized constants $C_1,C_2,C_3$:

\begin{eqnarray}\label{11.1}
\phi_1 = (3,1) \varphi_{2,1}^{var} -  (1,2) \varphi_{3,1}^{var} + C_1,
\end{eqnarray}
\begin{eqnarray}\label{11.2}
\phi_2 = (1,2) \varphi_{3,2}^{var} -  (2,3) \varphi_{1,2}^{var} + C_2,
\end{eqnarray}
\begin{eqnarray}\label{11.3}
\phi_3 = (2,3) \varphi_{1,3}^{var} -  (3,3) \varphi_{2,3}^{var} +  C_3.
\end{eqnarray}

The normalized constants $C_1,C_2,C_3$ are well-defined as follows. Take the constants 
$C_1,C_2,C_3$ in the formulas $(\ref{11.1})-(\ref{11.3})$ are equal to zero. 
This condition is equivalent to anyone of the two following conditions in the case
$(i+1,i) \ne 0$, $(i+2,i) \ne 0$, where the index $i$ is defined  $\pmod{3}$.
(In the case $(j,i)=0$, $j=i+1,i+2$, the corresponding condition is degenerate, but is valid.)
\begin{eqnarray}\label{phi}
\int_{U_i} \phi_i(\A^{0}_{i+1,i},\B_i) dU_i = 0, \quad \int_{U_i} \phi_i(\A^{0}_{i+2,i},\B_i) dU_i =0.
\end{eqnarray}

In the paper $\cite{A1}$ an alternative  gage of the constants $C_1,C_2,C_3$ in the formulas $(\ref{11.1})$---$(\ref{11.3})$
is presented. Namely, the normalized constants are defined from equations (it is sufficient to assume that the only equation, defined by the sum of the equations $(\ref{14})$-$(\ref{16})$ is satisfied): 
\begin{eqnarray}\label{14}
\int_{U_1} C_1 dU_1 = \int_{U_1} \varphi_{2,1}^{var} (\grad
\varphi_{3,1}^{var},\B_1) - (\grad \varphi_{2,1}^{var},\B_1)
\varphi_{3,1}^{var} dU_1 +
\end{eqnarray}
$$
\frac{2}{3} \int_{\R^3} \langle \A_1(x),\A_2(x),\A_3(x) \rangle d\R^3,
$$
\begin{eqnarray}\label{15} 
\int_{U_2} C_2 dU_2 = \int_{U_2}
\varphi_{3,2}^{var} (\grad \varphi_{1,2}^{var},\B_2) - (\grad
\varphi_{3,2}^{var},\B_2) \varphi_{1,2}^{var} dU_2 +
\end{eqnarray}
$$
\frac{2}{3} \int_{\R^3} \langle \A_1(x),\A_2(x),\A_3(x) \rangle d\R^3,
$$
\begin{eqnarray}\label{16} 
\int_{U_3} C_3 dU_3 = \int_{U_3}
\varphi_{1,3}^{var} (\grad \varphi_{2,3}^{var},\B_3) - (\grad
\varphi_{1,3}^{var},\B_3) \varphi_{2,3}^{var} dU_3 +
\end{eqnarray}
$$
\frac{2}{3} \int_{\R^3} \langle \A_1(x),\A_2(x),\A_3(x) \rangle d\R^3.
$$

The fraction $\frac{2}{3}$ of the second terms in the right side of the equations corresponds with the equation (10)\cite{A1},
in which a little simplification because of the new denotation  $(\ref{7})$ is presented.
This is proved in the following lemma.

\begin{lemma}\label{101}

--1. The equations $(\ref{14})$ -- $(\ref{16})$ determines the gauge corresponded to (10)\cite{A1}. 

--2. The equations $(\ref{14})$ -- $(\ref{16})$ depends no of a prescribed coordinate systems
$U_i = \Sigma_i \times S^1$ in the magnetic tubes.

\end{lemma}

\subsubsection*{Proof of Lemma $\ref{101}$}

Proof Statement 1. It is sufficient to proof, that the first term in the right side of the each equation
$(\ref{14})$-
$(\ref{16})$ is equal to the term $I_i-J_i$ in the formula  (13)\cite{A1}. Let us present the calculation in the case $i=1$, 
because the cases 
$i=2$, $i=3$ are obtained by the permutation of the indexes. Recall, that in \cite{A1} the following integrals are well defined: 
\begin{eqnarray}\label{I1}
I_1= \int(\B_1,[\A_2\varphi_{3,1} - \A_3\varphi_{2,1}])dU_1,
\end{eqnarray}
\begin{eqnarray}\label{J1}
J_1 = \int(\B_1,n)[\varphi_{2,1}lk_{3,1}fl_3 - \varphi_{3,1}lk_{2,1}fl_2]d\Sigma_1,
\end{eqnarray}
where $\varphi_{2,1}$, $\varphi_{3,1}$ are given by the equation $(\ref{7})$, $lk(3,1)$, $lk(2,1)$ are the linking coefficients 
of the central lines of the magnetic tubes, 
 $fl_2$, $fl_3$ are the flows trough the magnetic tubes $U_2$, $U_3$.
Obviously, the expressions
$(\ref{I1})$,  $(\ref{J1})$ depend of the constant $C_{j,i}$ in the formula  $(\ref{7})$,
but in calculations below this is not important.

Put  $\varphi_{2,1}$, $\varphi_{3,1}$ from the equations  $(\ref{7})$ to the equation $(\ref{I1})$, using the condition 
$(\ref{phi})$ for $i=1$. After simplifications we get: 
$$0= -I_1 + \int(\B_1,[\grad \varphi^{var}_{3,1}\varphi^0_{2,1} - \grad \varphi^{var}_{2,1}\varphi^{0}_{3,1}])dU_1 +$$
\begin{eqnarray}\label{cal2}
\int_{U_1} \varphi_{2,1}^{var} (\grad
\varphi_{3,1}^{var},\B_1) - (\grad \varphi_{2,1}^{var},\B_1)
\varphi_{3,1}^{var} dU_1.
\end{eqnarray}
The tird integral term in the right side of the formula coincides with the first term in the right side of the equation
 $(\ref{14})$. Let us prove that the first term in the second integral coincides with the first term
in $(\ref{J1})$. Using the Gauss-Ostrogradsky's formula for the ball $U_1 \setminus \Sigma_1$, we get that the difference 
of the values of the function 
 $\varphi^0_{3,1}$ restricted to the surface $-\Sigma_1$, and restricted to the surface  $\Sigma_1$, 
(each of the two surfaces is a subdomain in the ball $U_1 \setminus \Sigma_1$) is equal to  $lk(3,1)fl_3$. This proves the equality of the first two terms of the expressions 
$I_1$ and $J_1$. Analogous calculations for the firs two terms are satisfied. Statement 1 is proved.

Proof Statement 2. A formal proof follows from statement 1 and results $\cite{A1}$. Let us present an independent proof of the statement.
Assume that $j=1$ for simplifications of the denotations on the formula $(\ref{7})$. Then we get:
$$\varphi_{1,2}^0 \mapsto \varphi_{1,2}^0 + (1,2)\delta_{1},$$ 
$$\varphi_{1,3}^0 \mapsto \varphi_{1,3}^0 + (1,3)\delta_{1},$$ 
where $\delta_1$ is an arbitrary. Then for the function $\phi_1$ in $(\ref{11.1})$ we get the identity gauge:
\begin{eqnarray}\label{11.10}
\phi_1 \mapsto (3,1) \varphi_{2,1}^{var} + (3,1)(1,2) \delta_1 -  (1,2) \varphi_{3,1}^{var} - (3,1)(1,2)\delta_1 + C_1.
\end{eqnarray}
For the corresponding term in the right side of the formula $(\ref{14})$ we get:
\begin{eqnarray}\label{140}
\int_{U_1} \varphi_{2,1}^{var} (\grad
\varphi_{3,1}^{var},\B_1) - (\grad \varphi_{2,1}^{var},\B_1)
\varphi_{3,1}^{var} dU_1 \mapsto 
\end{eqnarray}
$$
\int_{U_1} \varphi_{2,1}^{var} (\grad
\varphi_{3,1}^{var},\B_1) - (\grad \varphi_{2,1}^{var},\B_1)
\varphi_{3,1}^{var} dU_1.
$$

Statement 2 and Lemma 
$\ref{101}$ is proved. 
\[   \]




Using the functions
 $\phi_i$ by the equations $(\ref{11.1})-(\ref{11.3})$, let us replace the vector-function  $\F$ in the formula 
$(\ref{FF})$ to a divergence-free vector function  $\F$ by the formula:
\begin{eqnarray}\label{FFF}
\F = (1,3)(2,3)\A_1 \times \A_2 + (2,1)(3,1)\A_2 \times \A_3 + (3,2)(1,2)\A_3 \times \A_1 
\end{eqnarray}
$$- \phi_1 \B_1(2,3) - \phi_2 \B_2(3,1) - \phi_3 \B_3(1,2). $$
Let us calculate the helicity of the vector-function
$(\ref{FFF})$, using the formula:
\begin{eqnarray}\label{MM}
 M(\B)= \int_{\R^3} (\G,\F) d\R^3 + \dots, 
\end{eqnarray}
where $\G$ is the vector-potential of  $\F$, $\dots$ in the formula $(\ref{MM})$ means the sum of the following  $10$ terms: 
\begin{eqnarray}\label{e}
e_{1,2,3}=-2(1,2)(2,3)(3,1)(\int_{\R^3} \langle \A_1,\A_2,\A_3 \rangle dx)^2,
\end{eqnarray}

\begin{eqnarray}\label{14.1}
f_{1} = -2(\int_{U_1} \varphi_{2,1}^{var} (\grad
\varphi_{3,1}^{var},\B_1) dU_1)(\int_{\R^3} \langle \A_1,\A_2,\A_3 \rangle dx), 
\end{eqnarray}
\begin{eqnarray}\label{15.1} 
f_2 = -2(\int_{U_2}
\varphi_{3,2}^{var} (\grad \varphi_{1,2}^{var},\B_2) dU_2)( \int_{\R^3} \langle \A_1,\A_2,\A_3 \rangle dx),
\end{eqnarray}
\begin{eqnarray}\label{16.1} 
f_3 = -2(\int_{U_3}
\varphi_{1,3}^{var} (\grad \varphi_{2,3}^{var},\B_3) dU_3 )( \int_{\R^3} \langle \A_1,\A_2,\A_3 \rangle dx),
\end{eqnarray}

\begin{eqnarray}\label{019}
d_{1,1} = -(2,3)^2 \int  \phi_1^2 (
\A_1,\B_1) dU_1,
\end{eqnarray}

\begin{eqnarray}\label{020}
d_{2,2} = -(3,1)^2 \int  
\phi_2^2(\A_2,\B_2) dU_2,
\end{eqnarray}

\begin{eqnarray}\label{021}
d_{3,3} = -(1,2)^2 \int  
\phi_3^2(\A_3,\B_3)  dU_3,
\end{eqnarray}

\begin{eqnarray}\label{022}
d_{1;3} = (2,3)(1,2) \int  
\phi_1^2(\A_3,\B_1) dU_1, 
\end{eqnarray}

\begin{eqnarray}\label{023}
d_{2;1} = (3,1)(2,3) \int  
\phi_2^2(\A_1,\B_2)  dU_2,
\end{eqnarray}

\begin{eqnarray}\label{024}
d_{3;2} = (1,2)(3,1) \int  
\phi_3^2(\A_2,\B_3) dU_3.
\end{eqnarray}

\begin{theorem}\label{A1}
The expression $(\ref{MM})$, in which the scalar potentials, given by the formulas  $(\ref{11.1})-(\ref{11.3})$ with the gauge conditions 
$(\ref{phi})$ is used, is the invariant of volume-preserved diffeomorphisms of the space. 
This invariant is not a function of pairwise integral integer coefficients of magnetic tubes. 
\end{theorem}

\subsubsection*{Proof of Theorem  $\ref{A1}$}
The statement of the theorem is a little modification of the main result of the paper \cite{A1}. The changes are following. Instead of 
the gauge
$(\ref{14})-(\ref{16})$ we use the gauge $(\ref{phi})$. This determines an extra 4 terms   $(\ref{e})$--$(\ref{16.1})$ 
in the formula of the invariant. Theorem  $\ref{A1}$ is proved.

\subsubsection*{Remark}
As far as the author knows  terms like $(\ref{019})$ - $(\ref{024})$ was initially adopted to higher helicity integrals by M.Berger in \cite{B}.  
\[  \]

The following lemma will be used in  the next section.

\begin{lemma} \label{srednee}
Let  $\theta_i \in S^1 = [0,2\pi)$, $i=1,2,3$,  be an angle, which parametrizes sections  $S_i(\theta)$ of the prescribed coordinate system in the magnetic tube
$U_i$, $i=1,2,3 \pmod{3}$, in particular, $S_i(0)=\Sigma_i$. Let us consider the gauge of the constants $C_i$ in the equations  $(\ref{11.1})$--$(\ref{11.1})$, which is given over the surface $S_i=S_i(\theta)$ by one of the following equations:
\begin{eqnarray}\label{Si1}
\int \phi_i(\A^{0}_{i+1,i},\B_i) dS_i = 0, \quad \int \phi_i(\A^{0}_{i+2,i},\B_i) dS_i =0.
\end{eqnarray}
(In the case $(i+1,i) \ne 0$, $(i+2,i)\ne 0$ the equations are equivalent. In the case
$(i+1,i) =0$, or
$i+2,i)=0$ the corresponding equation is the identity and can be omitted.) Assume that the invariant
 $M(\theta)$ is defined by the formula  $(\ref{MM})$, using the considered gauge. 
Then the mean value
$$\frac{1}{(2\pi)^3}\int_{0}^{2\pi} \int_0^{2\pi} \int_0^{2\pi} M(\theta_1,\theta_2,\theta_3)d\theta_1 d\theta_2 d \theta_3$$
 of the invariant $M(\theta_1,\theta_2, \theta_3)$, which is calculated using all prescribed  gages, which is parametrized over all possible sections $S_i(\theta)$, coincides with the integral $M$, which is calculated by the formula  $(\ref{MM})$ in the gauge $(\ref{phi})$.
\end{lemma}

\subsubsection*{Proof of Lemma  $\ref{srednee}$}

The integral $M$ is linear over the constants $C_1,C_2,C_3$. This property is satisfied because a changing of a gauge of the constants 
changes the integral $M$ by the first additive constant, which is proportional to the sum of the integrals  $(\ref{e})$ with the coefficient, and by  the second additive constant, which is proportional to the sum of the integrals $(\ref{14.1})$--$(\ref{16.1})$ with the coefficient $C_i$. 

Let $\phi_i$ be the scalar potential in the magnetic tube  $U_i$, which is defined using the gauge
$(\ref{phi})$. Then the function  
$\phi_i(\theta_i)$, which is defined using the gauge $(\ref{Si1})$, distinguishes from  $\phi_i$, which is defined using the gauge
$(\ref{phi})$, on the constant $\int \phi_i(\A^{0}_{i+1,i},\B_i) dS_i$. 
The mean value
$\theta_i$  of the family of the constant $C_i(\theta)$, by linearity, coincides with  the integral in the left side of the equation 
$(\ref{phi})$ and, therefore, is equal to zero. Lemma  $\ref{srednee}$ is proved.

\section{Invariant $M$ for a triple of closed magnetic lines} 

In this section we investigate the limit of the formula 
$(\ref{MM})$, if the thickness of the magnetic tubes tends to zero, the magnetic flow thought the transversal section of the each tube remains constant. In the limit we also assume that all magnetic lines inside the each tube are closed the magnetic flow over the time $t_0=2\pi$ inside each tube is the identity (this is a special case of Property 
$\ref{vosvr}$. The coordinate system  $U_i \cong S^1 \times D^2$ inside the each magnetic tube corresponds to the evolution map generated by the magnetic flow, end concentrates in the limit near the central line. Let us call a magnetic tube with this properties is elementary.

For an arbitrary 3 elementary tubes a 3-component link
$$\L = L_1 \cup L_2 \cup L_3 \subset \R^3, $$
generated by the central lines of the tubes,
is well-defined. Let us define by $\dot{x}_i$ the the tangent vector to the parametrized curve
$\L$ at the point $x_i
\in L_i \subset \L$, $i=1,2,3$. The magnetic line $L_i$ coincides with the limit position of elementary magnetic tubes.

For an arbitrary point
$x_i \in L_i$  define the vector-field 
$\A(x_i;x)$ with the singular point $x_i$ by the Bio-Savard formula: 

\begin{eqnarray}\label{BioSavard}
\A(x_i;x) = \frac{1}{4\pi} \frac{\dot{x}_i \times (x-x_i)}{(x-x_i)^2}.
\end{eqnarray}

Define the vector-field
$\A_i(x)$, called the vector-potential of the magnetic line $L_i$, which is singular over the curve  $L_i$, by the formula:

\begin{eqnarray}\label{Ai}
\A_i(x) = \oint_{L_i} \A(x_i;x)  dx_i, \quad x_i \in L_i.
\end{eqnarray}

For an arbitrary pair of points
$x_i \in L_i$, $x_j \in L_j$, $i \ne j$, let us define the vector-field
$\valpha_{i,j}(x_i,y_j;x)$ by the formula: 

\begin{eqnarray}\label{3}
\valpha(x_i,x_j;x) =  \A(x_i;x) \times \A(x_j;x).
\end{eqnarray}

Define the vector-field
$\valpha_{i,j}(x)$ with singularities on the lines  $L_i, L_j$ by the formula:

\begin{eqnarray}\label{4}
\valpha_{i,j}(x) = \A_i(x) \times
\A_j(x) = \oint_{L_i \cup L_j} \valpha_{i,j}(x_i,x_j;x) dx_i dx_j,
\quad x_i \in L_i, x_j \in L_j.
\end{eqnarray}


The vector-field
$(\ref{4})$ is equal to the limit of the potential  $\A_i \times \A_j$ in the formula $(\ref{FF})$ 
(or, in the formula \cite{A1},(19)). Denote by $(i,j)$, $i,j=1,2,3$, $i \ne j$, the integral linking coefficient of the $i$-th and $j$-th components  of $\L$. This coefficient is defined by the formula: 

\begin{eqnarray}\label{5}
(i,j) = \frac{1}{4\pi} \oint \oint \frac{\langle
\dot{x}_i,\dot{x}_j,x_i-x_j \rangle} {\|x_i-x_j\|^3} dx_i,
dx_j, \quad x_i \in L_i, x_j \in L_j,
\end{eqnarray}
because the magnetic flows over the transversal section of the magnetic line $L_i$ is normalized, the integral linking coefficient is an integer.

Denote by $\varphi_{j,i} : L_i \to \R^1$, $i \ne j$ a multivalued function, which is defined 
by the integration of 
$\A_j$ over the magnetic line $L_i$. This multivalued function is periodic with the period $(i,j)$ and determines by the formula:

\begin{eqnarray}\label{6}
\varphi_{j,i}(x_i) = \int_{pt_i}^{x_i} \A_j(x_i)  dx_i + C_{j,i}, \quad x_i
\in L_i,
\end{eqnarray}
where $pt_i \in L_i$ is the marked point on the line $L_i$.  By the construction 
we get $\varphi_{j,i}(pt_i)=C_{i,j}$ up to the additive constant  $(i,j)k$, $k \in \Z$. The multivalued function 
$\varphi_{j,i}$ is defined as the limit of the function  $(\ref{6})$ (see the formula (19)\cite{A1}, equation (13)\cite{A1}).

We have to define the function
 $\phi_i(x_i), \quad x_i \in L_i$, which is the limit of the corresponding function  $(\ref{11.1})$--$(\ref{11.3})$. 
Let us prove the following statement.

\begin{lemma}\label{cal}

Let $U_i$ be an elementary magnetic tube. 

--1. The equation  $(\ref{phi})$ implies the following: 
\begin{eqnarray}\label{U1}
\int_{U_i} \phi_i d U_i  = 0.
\end{eqnarray}

--2. The following equation is satisfied: 
\begin{eqnarray}\label{U1.2}
\int_{U_1} (\B_1,[\varphi_{2,1}^{var} \grad
\varphi_{3,1}^{var} - \varphi_{3,1}^{var}\grad \varphi_{2,1}^{var}]
) dU_1  =
\end{eqnarray}
$$ \int (\B_1,[\A_3\varphi_{2,1}-\A_2\varphi_{3,1}]) dU_1, $$
where a multivalued functions 
$\varphi_{2,1}$, $\varphi_{3,1}$ at the right side of the equation are given by the formula $(\ref{6})$,
for $C_{j,i}=0$, in particular, is equal to zero at the marked point  $pt \in L_i$. 
The right side of the equation is a single-valued function in 
$U_1$. The analogous equations for magnetic tubes  $U_2$ and $U_3$ are well-defined.

\end{lemma}

\subsubsection*{Proof of Lemma $\ref{cal}$}

Proof of Statement 1. In each the equation 
$(\ref{phi})$ (below for short the first equation is considered) 
the vectors $\A_{i+1}$, $\A_{i+2}$ point along the central line of the tube and corresponds to the prescribed parametrization, 
the function 
$(\A_{i+1}^0,\B_i)$ is the constant, the integral of this function over $U_i$  is equal to $(i+1,i)$.

Proof of Statement 2. 
Let us apply the equation  $(\ref{7})$ for an elementary magnetic tube. 
The function 
$\varphi_{j,i}^{0}(pt)$ is lifted to the linear function of $\R$-coordinate on the universal covering 
of the elementary magnetic tube, which is equal to zero at the marked point.  
The expression $(\ref{U1.2})$ is followed from 
$(\ref{cal2})$.
The right side of the equation  $(\ref{U1.2})$ is followed from the equation  $(\ref{I1})$.  The term  $(\ref{J1})$ is equal to zero, because
the condition 
$\varphi_{i+1,i}(pt)=\varphi_{i+2,i}(pt)=0$ implies analogous conditions on the cross-section surface $\Sigma_i$, which contains the marked point
$pt_i$. Statement 2 and Lemma  $\ref{cal}$ are proved.   
\[  \]

Denote by 
$p: \tilde L_i \to L_i$ the universal covering over the magnetic line $L_i$. The universal covering contains the marked point 
$\widetilde{pt} \in \tilde L_i \cong \R$, which is mapped into the marked point $pt_i \in L_i$. Let us consider a point $\tilde x_i \in \tilde L_i$,
denote the projection $p(\tilde x_i)$ by $x_i$, and consider the real $a_{j,i}=(\dot{x}_i,\A_j(x_i))$, where $\dot{x}_i$ is the tangent vector
to  $L_i$ at the point $x_i$, the vector-function
$\A_j$ is defined by the formula  $(\ref{Ai})$. Denote by $\theta_{j;x_i}: \tilde L_i \to \R$ the standard locally-constant function with
the jump at the point 
$x_i$: $\frac{\partial \theta_{j;x_i}}{\partial x}=a_{j,i}(j,i)$, $x_i \in \tilde L_i$, $i \ne j$, which is normalized as following:  $\theta_{j;x_i}(-\infty)=0$, where $(j,i)$ is the integral (integer) linking coefficient. 

Let us define the (formal) function
\begin{eqnarray}\label{tiltheta}
\varphi_{j,i} : \tilde L_i \to \R,
\end{eqnarray}
by the formula:
$$ \varphi_{j,i}(x_i) = \int_{x_i}^{+\infty} \theta_{j,x_i}dx. $$
Define the periodic function
\begin{eqnarray}\label{phi'}
\begin{array}{c}
\phi_i = \varphi_{i+1,i}(i+1,i)-\varphi_{i+2,i}(i+2,i),  \quad i=1,2,3;\quad
\lim_{\tilde x_i \to +\infty} \phi_i = 0.
\end{array}
\end{eqnarray}


Note that the integral $(\ref{tiltheta})$ is not convergent and the function $\varphi_{i,j}$ is not well-defined. The 
functions $(\ref{phi'})$ is well-defined and periodic, where the condition  $\lim_{\tilde x_i \to +\infty} \phi_i = 0$ is
assumed in the sense of its mean value. Using Lemma $\ref{srednee}$ the gauge $(\ref{phi'})$, which corresponds to $(\ref{U1})$, is calculated from the gage 
\begin{eqnarray}\label{phi2}
\phi_i(pt_i)=0, \quad  pt_i \in L_i
\end{eqnarray} 
by the variation of the marked point $pt_i$ over $L_i$.

The required function 
$(\ref{6})$ (the scalar potential) in the  gauge $(\ref{U1})$
is defined.
\[  \]



Let us rewrite the formula
$(\ref{MM})$ with the main term  $(\ref{FFF})$ and with the last terms  $(\ref{e})-(\ref{24})$ 
for magnetic fields inside 3 elementary magnetic tubes.

For an arbitrary 3 points
$x_i \in L_i, i=1,2,3$ let us define the vector-field $\F(x_1,x_2,x_3;x)$ 
with singular points
$x_1,x_2,x_3$ by the formula:

\begin{eqnarray}\label{17}
\F(x_1,x_2,x_3;x)=(2,3)(3,1) \valpha_{1,2}(x_1,x_2;x) + (3,1)(1,2)
\valpha_{2,3}(x_2,x_3;x) +
\end{eqnarray}
$$
(1,2)(2,3)\valpha_{3,1}(x_3,x_1;x),
$$
where $\valpha(x_i,x_j;x)$ is defined by the formula $(\ref{4})$.

Define the vector-field $\F(x)$ with singularities on the components of $\L$ by the formula: 

\begin{eqnarray}\label{17.1}
\F(x)= \oint_{\L} \F(x_1,x_2,x_3;x) dx_1 dx_2 dx_3.
\end{eqnarray}

The vector-field, which is given by the equation 
$(\ref{17.1})$, corresponds to the first 3 terms in the formula $(\ref{FFF})$.

Define the real number
 $W$ using the Gauss integral:  

\begin{eqnarray}\label{18}
W = \frac{1}{4 \pi} \int \int  \frac{ \langle \F(x), \F(y), (x-y)
\rangle}{ \|x-y \|^3} dx dy = \frac{1}{4\pi} \int \int \Gamma(\F(x),\F(y)) dx dy.
\end{eqnarray}
The integral $(\ref{18})$ arise when the integral $(\ref{MM})$ is calculated using the Gauss integral, 
this integral is decomposed into the sum of 6 simplest integrals, using the formula
$(\ref{17})$. The integral  $(\ref{18})$ is called the main term in the formula $M$ (see below $(\ref{40})$).

Let us define the vector field
$\A_i^{\phi}(x)$ with singularities on the magnetic line  $L_i$ by the formula: 

\begin{eqnarray}\label{Aphi}
\A^{\phi}_i(x) = \oint_{L_i} \phi_i(x_i)\A(x_i,x)  dx_i, \quad x_i
\in L_i.
\end{eqnarray}

Define the real
$b_{1;1,2}$, $b_{1;1,3}$,
$b_{2;2,3}$, $b_{2;2,1}$ $b_{3;3,1}$,  $b_{3;3,2}$ 
by the following integrals, using the field by the equations 
 $(\ref{Ai})$, $(\ref{Aphi})$:

\begin{eqnarray}\label{19}
b_{1;1,2} = -(2,3)^2(3,1) \int  \langle
\A_1(x),\A_2(x),\A_1^{\phi}(x) \rangle dx,
\end{eqnarray}

\begin{eqnarray}\label{20}
b_{1;1,3} = -(2,3)^2(1,2) \int  \langle
\A_3(x),\A_1(x),\A_1^{\phi}(x) \rangle dx,
\end{eqnarray}

\begin{eqnarray}\label{21}
b_{2;2,3} = -(3,1)^2(1,2) \int  \langle
\A_2(x),\A_3(x),\A_2^{\phi}(x) \rangle dx,
\end{eqnarray}

\begin{eqnarray}\label{22}
b_{2;2,1} = -(3,1)^2(2,3) \int  \langle
\A_1(x),\A_2(x),\A_2^{\phi}(x) \rangle dx,
\end{eqnarray}

\begin{eqnarray}\label{23}
b_{3;3,1} = -(1,2)^2(2,3) \int  \langle
\A_3(x),\A_1(x),\A_3^{\phi}(x) \rangle dx,
\end{eqnarray}

\begin{eqnarray}\label{24}
b_{3;3,2} = -(1,2)^2(3,1) \int  \langle
\A_2(x),\A_3(x),\A_3^{\phi}(x) \rangle dx.
\end{eqnarray}

Let us remark, that the integrals
$(\ref{19})$--$(\ref{24})$, generally speaking, are not trivial, because the vector fields
$\A_i$ and $\A_{i}^{\phi}$, are linear independent. The sum of the integrals
$(\ref{19})$--$(\ref{24})$ corresponds to the second part of the main term in the integral $(\ref{MM})$, which is calculated using the Gauss integral.
The terms are calculated using the Gauss integral for the pair of the vector fields, the first vector field in the pair is the one of the term in
$(\ref{FFF})$ (for example, the term  $(1,2)(2,3)\A_1 \times \A_3$), the second vector field in the pair is given by one of the 3 last terms in 
$(\ref{FFF})$ (for example, by the term  $\phi_1 \B_1$ or $\phi_3 \B_3$).

Let us define the reals
$b_{1;2,3}$, $b_{2;3,1}$,
$b_{3;1,2}$ by the following integrals: 

\begin{eqnarray}\label{25}
b_{1;2,3} = -(1,2)(2,3)(3,1) \int  \langle
\A_2(x),\A_3(x),\A_1^{\phi}(x) \rangle dx,
\end{eqnarray}

\begin{eqnarray}\label{26}
b_{2;3,1} = -(1,2)(2,3)(3,1) \int  \langle
\A_3(x),\A_1(x),\A_2^{\phi}(x) \rangle dx,
\end{eqnarray}

\begin{eqnarray}\label{27}
b_{1;2,3} = -(1,2)(2,3)(3,1) \int  \langle
\A_1(x),\A_2(x),\A_3^{\phi}(x) \rangle dx.
\end{eqnarray}

The sum of the integrals
$(\ref{25})$--$(\ref{27})$ is defined using the Gauss integral to calculate the main term in the formula 
$(\ref{MM})$ for the corresponding pair of vector fields. 
As in the cases 
$(\ref{19})$--$(\ref{24})$, the first term in the pair coincides with one of the first 3 terms in  $(\ref{FFF})$ (
for example, the term
$(1,2)(2,3)\A_1 \times \A_3$), the second term is a one of the last of 3 terms in the formula
$(\ref{FFF})$ (for example, the term $\phi_2 \B_2$).

Define the reals
 $c_{1;1}$, $c_{2;2}$, $c_{3;3}$,
$c_{1;2}$, $c_{2;3}$, $c_{3;1}$ by the following integrals: 

\begin{eqnarray}\label{28}
c_{1;1} = (2,3)^2 \oint \phi_1  (\dot{x}_1,\A^{\phi}_1)  dx_1,
\end{eqnarray}

\begin{eqnarray}\label{29}
c_{2;2} = (3,1)^2 \oint \phi_2  (\dot{x}_2,\A^{\phi}_2)  dx_2,
\end{eqnarray}

\begin{eqnarray}\label{30}
c_{3;3} = (3,1)^2 \oint \phi_3  (\dot{x}_3,\A^{\phi}_3)  dx_3,
\end{eqnarray}

\begin{eqnarray}\label{31}
c_{1;2} = 2(2,3)(3,1) \oint \phi_2  (\dot{x}_2,\A^{\phi}_1)
dx_2,
\end{eqnarray}

\begin{eqnarray}\label{32}
c_{2;3} = 2(3,1)(1,2) \oint \phi_3  (\dot{x}_3,\A^{\phi}_2)
dx_3,
\end{eqnarray}

\begin{eqnarray}\label{33}
c_{3;1} = 2(1,2)(2,3) \oint \phi_1  (\dot{x}_1,\A^{\phi}_3)
dx_1.
\end{eqnarray}

The integrals
$(\ref{25})$--$(\ref{27})$ are also obtained using Gauss integral in the formula $(\ref{MM})$ for a pair of vector fields. 
Both vector-fields in the pair are defined are defined by last terms in $(\ref{FFF})$. 

Define the reals
$d_{1;1}$, $d_{2;2}$, $d_{3;3}$,
$d_{1;2}$, $d_{2;3}$, $d_{3;1}$ by the following integrals: 

\begin{eqnarray}\label{34}
d_{1;1} = -(2,3)^2 \oint \phi^{2}_1  (\dot{x}_1,\A_1) dx_1,
\end{eqnarray}

\begin{eqnarray}\label{35}
d_{2;2} = -(3,1)^2 \oint \phi^{2}_2  (\dot{x}_2,\A_2) dx_2,
\end{eqnarray}

\begin{eqnarray}\label{36}
d_{3;3} = -(3,1)^2 \oint \phi^{2}_3  (\dot{x}_3,\A_3) dx_3,
\end{eqnarray}

\begin{eqnarray}\label{37}
d_{1;2} = (3,1)(2,3) \oint \phi^{2}_2  (\dot{x}_2,\A_1) dx_2,
\end{eqnarray}

\begin{eqnarray}\label{38}
d_{2;3} = (1,2)(3,1) \oint \phi^{2}_3  (\dot{x}_3,\A_2) dx_3,
\end{eqnarray}

\begin{eqnarray}\label{39}
d_{3;1} = (2,3)(1,2) \oint \phi^{2}_1  (\dot{x}_1,\A_3) dx_1.
\end{eqnarray}
The integrals
$(\ref{24})$--$(\ref{39})$ corresponds to last terms  $(\ref{19})$--$(\ref{24})$  in the formula  $(\ref{MM})$.
Not difficult to check that $c_{i,i} = - d_{i,i},\quad i=1,2,3$.

Define the terms, which correspond to 
$(\ref{14.1})$--$(\ref{16.1})$, by the equations: 

\begin{eqnarray}\label{14.2}
f_{1} = -2(\oint_{L_1} (\dot{x}_1,[\A_3\varphi_{2,1}-\A_2\varphi_{3,1}]) dL_1)(\int_{\R^3} \langle \A_1,\A_2,\A_3 \rangle dx), 
\end{eqnarray}
\begin{eqnarray}\label{15.2} 
f_2 = -2(\oint_{L_2}(\dot{x}_2,[\A_1\varphi_{3,2}-\A_3\varphi_{1,2}])
 dL_2)( \int_{\R^3} \langle \A_1,\A_2,\A_3 \rangle dx),
\end{eqnarray}
\begin{eqnarray}\label{16.2} 
f_3 = -2(\oint_{L_3}(\dot{x}_3,[\A_2\varphi_{1,3}-\A_1\varphi_{2,1}])dL_3)
( \int_{\R^3} \langle \A_1,\A_2,\A_3 \rangle dx).
\end{eqnarray}
The correspondence with $(\ref{14.1})$--$(\ref{16.1})$ is satisfied by Lemma $\ref{cal}$, statement 2.
To simplify the integration over the variation of marked points $pt_i$ we may assume that in the formula $(\ref{14.1})$--$(\ref{16.1})$
the gage for $\varphi_{i,j}$ is given by the formal asymptotic $(\ref{tiltheta})$, this gives the gauge for $[\A_3\varphi_{2,1}-\A_2\varphi_{3,1}]$
and the analogous two terms.

Define the term, which corresponds to 
 $(\ref{e})$, by the formula:  
\begin{eqnarray}\label{71}
e(1;2;3) = -2(1,2)(2,3)(3,1)(\int_{\R^3} \langle \A_1(x),\A_2(x),\A_3(x) \rangle dx)^2.
\end{eqnarray}
Note that vector-potentials  $\A_i(x)$ in $(\ref{71})$ are defined by the formula $(\ref{Ai})$.

\subsubsection*{Formula of the invariant $M$ for elementary magnetic tubes} 


Define the real
 $M=M(\L)$ by the following sum of the integrals: 
\begin{eqnarray}\label{40}
M = W + \sum_{i;j=1}^3 b_{i;i,j} + \sum_{i=1}^3 (b_{i;i+1,i+2} +
c_{i,i+1} + d_{i;i+1} + f_i) + e(1,2,3).
\end{eqnarray}

In the formula
$(\ref{40})$  a parametrization  of  $\L$ is given by a circle of an arbitrary radius.

\begin{theorem}$\label{th1}$

1. The terms in the right side of $(\ref{40})$ is absolutely integrable.

2. The integral $(\ref{40})$  defines an invariant of isotopy classes of 
$\L$. The invariant is not a function of pairwise integral linking numbers of components. 
\end{theorem}

\subsubsection*{Proof of Theorem  $\ref{th1}$}
Proof of Statement 1.   Let us present an explicit estimation of the main term $W$, which is given by the equation $(\ref{18})$. This term is divided into 2 parts. 
The first group  contains the term:
\begin{eqnarray}\label{80}
(2, 3)^2(3, 1)^2 \oint_{L_1} \oint_{L'_1} \oint_{L_2} \oint_{L'_2} \int_{\R^3 \times \R^3}
\end{eqnarray}
$$
\frac{\langle \A(x_1; x) \times \A(x_2 ; x),
\A(x'_1; x') \times \A(x'_2 ; x'),(x-x')\rangle}{\| x-x'\|^3} dx dx' dx_1 dx'_1 dx_2 dx'_2,
$$
and analogous two terms, which are obtained by the cyclic permutation of the indexes.

The second group  contains the term:
\begin{eqnarray}\label{81}
2(3, 1)(1, 2)(2,3)^2 \oint_{L_1} \oint_{L'_1} \oint_{L_2} \oint_{L_3} \int_{\R^3 \times \R^3}
\end{eqnarray}
$$
\frac{\langle \A(x_1; x) \times \A(x_2 ; x),
\A(x'_3; x') \times \A(x'_1 ; x'),(x-x')\rangle}{\| x-x'\|^3}dx_1 dx'_1 dx_2 dx'_3,
$$
and  analogous two terms, which are obtained by the cyclic permutation of the indexes.

In the expression
$(\ref{80})$ the first term  $(2, 3)^2(3, 1)^2$ is the order 4 polynomial of pairwise linking numbers of the components.
The main term
$(\ref{80})$ when 
$x_1 \to x'_1$ is regular: the singularity is of the order $r^{-3}$ over $\R^4$.  
Estimations of the term
$(\ref{81})$ and of analogous terms are similar. This proves that the main term $W$ is convergent. 

The term $b_{1,1;2}$ and the analogous two terms in the formula   $(\ref{40})$
are given by the integral, which has the formal pole of the order $r^{-4}$ over 4-dimensional space.
The integral is converge, because the order of the pole is $1$  less then its formal order.

 The convergence of other terms in the formula
  $(\ref{40})$ is evident. Statement 1 of Theorem $\ref{th1}$ is proved.
  
Proof of Statement 2. 
Let us consider one-parameter family of the magnetic fields
 $\B(\varepsilon)$ over the parameter $\varepsilon \to 0+$, in this family  thicknesses of  tubes
 tend to zero and the integral flow over the cross-sections of  tubes is preserved.
 
The integral $(\ref{MM})$ is well-defined  over $\R^3$ and the limit
\begin{eqnarray}\label{3000}
\lim_{\varepsilon \to 0+} M(\B(\varepsilon))
\end{eqnarray}
is well-defined.

The expression $(\ref{40})$ is defined as the result of the limit  $(\ref{3000})$.
For an arbitrary marked points on the magnetic lines the gage $(\ref{phi'})$ (but not the gage $(\ref{U1})$)
is satisfied. The relation between those two gage  is given by Lemma 
$\ref{srednee}$.

In the limit  $(\ref{3000})$ each term in 
$(\ref{MM})$ is regular, in particular, the invariant $M$ is not degenerated. 
By Theorem 
$\ref{A1}$ we get the proof of Statement 2 of Theorem $\ref{th1}$.

Theorem $\ref{th1}$ is proved.
\[  \]

\subsubsection*{Remark}
The invariant
 $M$ is skew-symmetric with respect to the mirror symmetry of the space. The dimension of $M$ is equal to $G^{12}sm^{6}$.
\[  \]

\section{Invariant $M$ is a Vassiliev's invariant}

In  this section we deduces from Theorem $\ref{th1}$ the following corollary.

\begin{theorem}\label{Vassinv}
The invariant $M$, see the formula $(\ref{40})$ is a finite-type invariant in the sense of V.A.Vassiliev.
\end{theorem}

Let $\L \subset \R^3$ be a 3-component link, let $Conf^{r}(\L) = (\L)^r$ be a finite-type (non-connected) configuration space, associated with the Cartesian product
of $r$ copies of $\L$ (comp. with definition of configuration spaces in \cite{C-K-Sh}).

Let 
\begin{eqnarray}\label{conf}
\varphi: Conf^r(\L) \to \R 
\end{eqnarray}
be an integrable function. 

\begin{lemma}\label{Vass}
All terms in the right side of the equation $(\ref{40})$ is defined as an integral of a suitable function $(\ref{conf})$.
\end{lemma}

\subsubsection*{Proof of Lemma $\ref{Vass}$}
For the terms $(\ref{18})$, $(\ref{71})$ Lemma is evident.  Let us prove the statement for the terms $(\ref{19})$, $(\ref{31})$,
$(\ref{37})$, and $(\ref{14.2})$. 

The integral linking numbers $1,2)$ $(2,3)$, $(3,1)$ 
are, evidently, given by the corresponding integrals over the configuration spaces, and we may assume that 
$(2,3)$, $(3,1)$ in the formula $(\ref{19})$ and $(1,2)$, $(3,1)$ in the formula $(\ref{phi'})$ are constant. 
The  integral
\begin{eqnarray}\label{190}
\int  \langle
\A_1(x),\A_2(x),\A_1^{\phi}(x) \rangle dx,
\end{eqnarray}
is well-defined over the subspace $Conf^7 = \L^7 \supset L_1^4 \times L_2^2 \times L_3$ as following.
Take two points $x_1, x_2$ on $L_1$ and associate the vectors fields $\A_1(x_1;x)$, $\A_1(x_2;x)$ as the first and the third vector of the triple
in the formula $(\ref{190})$.
Take the oriented segment $[x_1,x_2]$ from $x_1$ toward $x_2$ on $L_1$ and take two points $x_3, x_4 \in [x_1,x_2]$.
Take the points $y_1 \in L_2$, $z_1 \in L_3$, associate the vectors $\A_2(y_1;x_3)$, $\A_3(z_1;x_4)$  to calculate the integrals
$\varphi_2$, $\varphi_3$   over $[x_1,x_2]$ as is required in the formula
$(\ref{phi'})$. Take the last point $y_2 \in L_2$ to calculate the second vector $\A_2(y_2;x)$ of the triple in the formula $(\ref{190})$. 
Lemma for the term $(\ref{19})$ is proved.

The integral
\begin{eqnarray}\label{310}
\oint \phi_2  (\dot{y}_2,\A^{\phi}_1)
dy_2, \quad y_2 \in L_2
\end{eqnarray}
is well-defined over the subspace $Conf^{11} = \L^{11} \supset L_1^5 \times L_2^4 \times L_3^2$ as following.
Take marked points $pt_1 \in L_1$, $pt_2 \in L_2$, these points are variated. The vector $\A^{\phi}_1$ is well-defined as following: take 
the point $x_1 \in L_1$ and associated with this point the vector-field $\A_1(x_1;x)$, and take points $x_2$, $x_3$ on the segment 
$[pt_1,x_1]$ and points $y_1 \in L_2$, $z_1 \in L_3$ to define the function $\phi_1$ and the vector $\A^{\phi}_1$ by the integration
over $[pt_1,x_1]$ as in the case  for the term $(\ref{19})$ above. Take a point $y_2 \in L_2$. Take points $y_3, y_4 \in [pt_2,y_2]$, $x_4 \in L_1$, $z_2 \in L_3$ to calculate $\phi_2$.  Calculate the function $\phi_2  (\dot{y}_2,\A^{\phi}_1)$ at the point $y_2 \in L_2$ in the formula $(\ref{190})$. 
Lemma for the terms $(\ref{31})$ is proved.

For the terms $(\ref{37})$, $(\ref{14.2})$ the proof is analogous. Lemma $\ref{Vass}$ is proved.

\subsubsection*{Proof of Theorem $\ref{Vassinv}$}
The proof follows from Lemma $\ref{Vass}$. 
Let $\varepsilon > 0$ be an arbitrary.  Let us prove that the algebraic sum of values of $M$ invariant over a collection $\aleph=\aleph(\varepsilon)$, of  $2^k$ links, which is the resolution a singular link $\L$ with $k$ self-crossings $\{a_i\}$, $i=1, \dots, r$, equals to zero. 

Assume that the number $k$ of self-crossings of a singular link $\L$ is greater then the number $r$ of points in the configuration space.
Then there exists a finite collection of neighborhoods $U_i \subset Conf^r(\L)$, $i=1, \dots, i_{max}$, such that $\cup_i U_i = Conf^r(\L)$ and for an arbitrary
$U_i$ 
there exists a  self-crossing point $a$ of the collection $\{a_i\}$, such that the regular $\varepsilon$--neighborhood of $a$ contains no point of $U_i$. Therefore
if the support of $\aleph(\varepsilon)$ is sufficiently small, the integral contribution of points from $U_x$ into the algebraic sum $\sum_{\alpha \in \aleph}M(\alpha))$ is trivial. Therefore we get $\sum_{\alpha \in \aleph}M(\alpha)) =0$. Theorem $\ref{Vassinv}$ is proved.

\section{Invariant $M$ is an asymptotic finite-order functional}

To prove Main Theorem
\label{mainTh} the following result is required. 

\begin{lemma}\label{776}

1. The main term  $(\ref{18})$ in the expression of $M$ is an asymptotic functional of magnetic fields
on the space 
 $\Omega$ in the sense of Definition $\ref{local}$.


2. The term $(\ref{37})$ the analogous property is satisfied.
\end{lemma}

\subsubsection*{Proof of Lemma $\ref{776}$}


Prof of Statement 1. 
Let us start with the definition of the term 
$(\ref{80})$ for a generic magnetic field  $\B \in \Omega$ using configuration space, for which the condition of Theorem $\ref{th4}$  is satisfied.
The coordinates of a point of the space
$K_{3,4;2}$ are given by collections $\{l_1,t_{1,1}, \dots t_{1,4},l_{2},t_{2,1}, \dots t_{2,4},l_3,t_{3,1}, \dots t_{3,4};y_1,y_2 \}$,
where $l_i \in U_i$, $t_{i,j} \in [0,T] \subset \R_{i,j}$, $j=1,2,3,4$, $y_1,y_2 \in \R^3$. 

The evaluation map
$F: K_{3,4;2} \to U_1^{4} \times U_2^4 \times U_3^4 \times (\R^3)^2$ is defined by the formula
 $$ F(l_1,t_{1,1}, \dots t_{1,4},l_{2},t_{2,1}, \dots t_{2,4},l_3,t_{3,1}, \dots t_{3,4};y_1,y_2 ) = $$
 $$(g^{t_{1,1}}(l_1), \dots g^{t_{1,4}}(l_1), g^{t_{2,1}}(l_2), \dots g^{t_{2,4}}(l_2), g^{t_{3,1}}(l_3), \dots g^{t_{3,4}}(l_3)),y_1,y_2$$
where $g^t$ is the magnetic flow. Therefore, the space 
$K_{3,4;2}$ is the configuration space of ordered collections of 4- points 
on the first magnetic line, issued from $l_1$, of 4 points on the second magnetic line, issued from $l_2$, and of 4 points on the third magnetic line, issued from $l_3$. The space is extended by an ordered pair of points
$(y_1,y_2) \in (\R^3)^2$. On the configuration space
$K_{3,4;2}$ the standard volume form $dK_{3,4;2}$ is well-defined.

To apply Theorem $\ref{th4}$ let us extend the configuration space $K_{3,4;2}$ by the following Cartesian product:
\begin{eqnarray}\label{Cart}
K^{1,2,3}_{3,2+1+1} \times K^{1,2}_{2,1} \times K^{2,3}_{2,1} \times \bar{K}^{2,3}_{2,1} \times K^{3,1}_{2,1}, 
\end{eqnarray}
where $K_{3,2+1+1}^{1,2,3} = \{l_1,t_{1,1},t_{1,2},l_{2},t_{2,1},l_3,t_{3,1};y_1,y_2 \}$, $K^{1,2}_{2,1}= \{l_1,t_{1,3},l_{2},t_{2,2}\}$,
$K^{2,3}_{2,1}= \{l_2,t_{2,3},l_{3},t_{3,3}\}$, $\bar{K}^{2,3}_{2,1}= \{l_2,t_{2,4},l_{3},t_{3,4}\}$, $K^{3,1}_{2,1}= \{l_1,t_{1,4},l_{3},t_{3,2}\}$. 
The tautologous projection 
\begin{eqnarray}\label{Cartproj}
K^{1,2,3}_{3,2+1+1} \times K^{1,2}_{2,1} \times K^{2,3}_{2,1} \times \bar{K}^{2,3}_{2,1} \times K^{3,1}_{2,1} \to K_{3,4;2}
\end{eqnarray}
is well defined.

Let us define the function 
\begin{eqnarray}\label{W342}
W_{3,4;2}: K_{3,4;2} \to \R, 
\end{eqnarray}
which is called the density of the functional $W$ for the corresponding summand in the formula $(\ref{18})$
(this function is a specification of the function $(\ref{Fq})$). The function $W_{3,4;2}$ is lifted by the projection $(\ref{Cartproj})$
to the product of the following functions:
\begin{eqnarray}\label{W1}
W^{1,2,3}: K_{3,2+1+1}^{1,2,3} \to \R,
\end{eqnarray}
\begin{eqnarray}\label{W2}
\begin{array}{c}
W^{1,2}: K^{1,2}_{2,1} \to \R, \quad W^{2,3}: K^{2,3}_{2,1} \to \R, \\ \bar{W}^{2,3}: K^{2,3}_{2,1} \to \R, \quad W^{3,1}: K^{3,1}_{2,1} \to \R.
\end{array}
\end{eqnarray}

Let us define the function $(\ref{W1})$. With a quadruple of points 
$x_{1,1}=g^{t_{1,1}}(l_1), x_{1,2}=g^{t_{1,2}}(l_1), x_{2,1}=g^{t_{2,1}}(l_2), x_{3,1}=g^{t_{3,1}}(l_3)$ 
the following integral kernel is associated:
\begin{eqnarray}\label{jadro}
W^{1,2,3}(x_{1,1},x_{2,1},x_{1,2},x_{3,1};y_1,y_2)=\Gamma(\valpha_{1,2}(x_{1,1},x_{2,1}),\valpha_{1,3}(x_{1,2},x_{3,1});y_1,y_2),
\end{eqnarray}  
where $\Gamma(y_1,y_2)$ is the integral kernel in the Gauss integral $(\ref{18})$, which is given by
\begin{eqnarray}\label{18}
\Gamma(y_1,y_2) = \frac{1}{4 \pi} \frac{ \langle \valpha_{1,2}(x_{1,1},x_{2,1}),\valpha_{1,3}(x_{1,2},x_{3,1}) , (y_1-y_2)
\rangle}{ \|y_1-y_2 \|^3},
\end{eqnarray}
the vector-field $\valpha_{1,2}(x_{1,1},x_{2,1})$ with singular points $x_{1,1},x_{2,1}$ is given by the expression $(\ref{3})$ for $x_i=x_{1,1}, x_j=x_{2,1}$, the vector-field $\valpha_{1,3}(x_{1,2},x_{3,1})$ with singular points $x_{1,2},x_{3,1}$ is given by the expression $(\ref{3})$ for $x_i=x_{1,2}, x_j=x_{3,1}$.

The function $(\ref{W342})$ is polynomial of the following collection of function: 

--1 coordinates of the vector $\valpha_{1,2}(x_{1,1},x_{2,1})$,

--2 coordinates of the vector $\valpha_{1,3}(x_{1,2},x_{3,1})$,

--3  components of the correlation tensor, which is defined by the kernel in the Gauss integral $(\ref{18})$ (see $(\ref{jadro})$), 
where the arguments are given by the base vectors at the points $y_1$, $y_2$. 

Each function, described by the list 1-3 is integrable (belongs to $L^1(\R^3)$).
This proves that the function $(\ref{W342})$ satisfies  Condition 1 of Theorem $\ref{th4}$.

Let us prove Condition 2 of Theorem $\ref{th4}$ for the function $(\ref{W342})$.
Define the function:
\begin{eqnarray}\label{intjadro}
A(x_{1,1},x_{1,2},x_{2,1},x_{3,1})=\int \int W^{1,2,3}(x_{1,1},x_{1,2},x_{2,1},x_{3,1};y_1,y_2) dy_1 dy_2,
\end{eqnarray}  
see $(\ref{Fq})$.
Let us construct 
a limiting tensor $a^{[2]}: (U_1)^2 \times U_2 \times U_3 \to \R_+$, for the function $(\ref{intjadro})$, see $(\ref{a})$. 
Let us define a real parameter  $\delta(x_{1,1},x_{1,2},x_{2,1},x_{3,1})$ as the smallest distance from two short $\varepsilon$-segments $I_{1,1}$, $I_{1,2}$ on the magnetic line  $L_1$, $x_{1,1} \in I_{1,1}$, $x_{1,2} \in I_{1,2}$, and two short $\varepsilon$-segments $I_2 \subset L_2$, $x_{2,1} \in I_2$,  $I_3 \subset L_3$, $x_{3,1} \in I_3$.

Analogously to the construction of \cite{A3} Section 3 let us define the $\varepsilon$-smoothing of the function $(\ref{intjadro})$ by the following integral:
\begin{eqnarray}\label{intjadroeps}
A_{\varepsilon}(x_{1,1},x_{1,2},x_{2,1},x_{3,1})=\int_{-\varepsilon}^{+\varepsilon} \int_{-\varepsilon}^{+\varepsilon} \int_{-\varepsilon}^{+\varepsilon} \int_{-\varepsilon}^{+\varepsilon} A(z_{1,1},z_{1,2},z_{2,1},z_{3,1}) dz_{1,1} dz_{1,2} dz_{2,1} dz_{3,1},
\end{eqnarray} 
where $x_{1,1}-\varepsilon < z_{1,1} < x_{1,1}+\varepsilon$, $x_{1,2}-\varepsilon < z_{1,2} < x_{1,2}+\varepsilon$,
$x_{3,1}-\varepsilon < z_{3,1} < x_{3,1}+\varepsilon$, $x_{3,1}-\varepsilon < z_{3,1} < x_{3,1}+\varepsilon$.
 
\begin{lemma}$\label{77}$

There exist $\varepsilon_0 > 0$, $\delta_0 > 0$, $C>0$  which are depended only on the $C^2$--norm of $\B$ in $U$,
such that there exists $\varepsilon > 0$, (without loss of a generality one may assume that $\varepsilon_0 << \varepsilon$, and $\delta_0 << \varepsilon$), such that for arbitrary 2 points $x_{1,1} \in L_1, x_{1,2} \in L_1$, $\dist(x_{1,1};x_{1,2})<\varepsilon_0$, 
and for arbitrary 2 points $x_{2,1} \in L_2$, $x_{3,1} \in L_3$, and $\delta(x_{1,1},x_{1,2},x_{2,1},x_{3,1})<\delta_0$ the function $(\ref{intjadro})$ satisfies following inequality:
\begin{eqnarray}\label{avarepsiloneps}
A_{\varepsilon}(x_{1,1},x_{1,2},x_{2,1},x_{3,1}) <  C \ln(\delta^{-1}).
\end{eqnarray}

2. There exist $\varepsilon_0 > 0$, $\delta_0 > 0$, $C>0$, which are depended only on the $C^2$--norm of $\B$ in $U$,
such that there exists $\varepsilon > 0$, $\varepsilon = \varepsilon(\varepsilon')$ (without loss of a generality one may assume that $\delta_0 << \varepsilon$), such that for arbitrary 4 points $x_{1,1} \in L_1, x_{1,2} \in L_1$, $x_{2} \in L_2, x_{3} \in L_3$ with the maximal distance between the points
not less then $\varepsilon_0$, and $\delta(x_{1,1},x_{1,2},x_{2,1},x_{3,1})<\delta_0$  the function $(\ref{intjadro})$ satisfies the inequality $(\ref{avarepsiloneps})$.

3. There exist $\varepsilon_0 > 0$, $C>0$, which are depended only on the $C^2$--norm of $\B$ in $U$,
there exists $\varepsilon > 0$, $\varepsilon = \varepsilon(\varepsilon')$, such that for arbitrary 4 points $x_{1,1} \in L_1, x_{1,2} \in L_1$, $x_{2} \in L_2, x_{3} \in L_3$ with the minimal distance between the points
not less then $\varepsilon_0$ the function $(\ref{intjadro})$ satisfies the inequality:
\begin{eqnarray}\label{avarepsilonepss}
A_{\varepsilon}(x_{1,1},x_{1,2},x_{2,1},x_{3,1}) <  C.
\end{eqnarray}

\end{lemma}

\subsubsection*{Proof of Lemma $\ref{77}$}

Proof of Statement 1.
Assume with no loss of a generality that the magnetic field
in $U' \subset U$, $I_{1,1},I_{1,2},I_2,I_3 \subset U'$, is represented by a  vector of the unite light and that $\varepsilon=1$. Define a dimensionless parameter
 $\lambda \in [1,+\infty)]$, which is the fraction of the distance between the magnetic lines  $(L_2,L_3)$ to 
 the distance between the magnetic lines  $(L_2,L_1)$, assuming that the order of the line corresponds to $\lambda \ge 1$.
denote the distance between the segments
$I_1$ and $I_2$  by $\delta$, then the distance between $L_2$ and $L_3$ not less then  $\delta$. 

Let us investigate the integral  $(\ref{intjadroeps})$
for the given configuration of segments of magnetic lines, which is defined by the parameters 
$\delta$.
Let us zoom-out the scale in 
 $\delta^{-1}$ times, the length of the each magnetic line will be equal to $\delta^{-1}$, $\delta^{-1} >> 1$. 
With respect to this transformation the each  elementary  potential is transformed (with a switch-out of a density) by the coefficient
$\delta^3$. At the end configuration the each elementary dipoles remains of the order 1. 

Let us show that the absolute value of the integral  
 $(\ref{intjadroeps})$  is estimated from below 
over the parameter 
  $\delta \to 0+$ by a value of the order $\delta^{-1}\ln{\delta^{-1}}$. The estimation 
is given by the product of

-- $\delta^{8}$, the dimension of the $4$ elementary dipoles,

-- $\delta^{-6}$, the dimension of the volume form, which is required to calculate the integral kernel in $(\ref{jadro})$, 

--$\delta^{-4}$ the dimension of the length forms of 4 segments of magnetic lines,

--$\delta^{2}$ the dimension of  $\vert y_1-y_2 \vert^{-2}$, which is required to calculate the Gauss integral.

Denote by $\B$ the original magnetic field, associated with the $\varepsilon$-segments $I_{1,1},I_{1,2},I_2,I_3$ of the magnetic lines $L_1,L_2,L_3$ and by $\B'$ the magnetic field after the transformation, associated with the $\varepsilon\delta^{-1}$--segments $I'_{1,1},I'_{1,2},I'_2,I'_3$
of the magnetic lines $L_1',L_2',L_3'$. Denote by $A'_{\varepsilon}$ result of the transformation of the expression $(\ref{intjadroeps})$.
By the calculation above we get $A'_{\varepsilon} = A_{\varepsilon}$. By the straightforward calculation we get $A'_{\varepsilon}$ is of the order $\delta^{-1}\ln{\delta^{-1}}$, if $\B$ and $\B'$ are of the order $1$. 
But the formal asymptotic is one order more then the value of the integral $A_{\varepsilon}$, because if all the $4$ vectors of $\B$ are parallel to
the common vertical axis we get $A_{\varepsilon}=0$ by the mirror symmetry. Statement 1 is proved.

Proof of Statement 2. Take $\varepsilon_0$ is much less then $\varepsilon_0$,$\delta_0$, which are defined in  Statement 1.
By analogous arguments we get by a formal asymptotic  the required inequality $(\ref{intjadroeps})$. Statement 2 is proved.

Statement 3 is evident. 
Lemma $\ref{77}$ is proved.

\subsubsection*{Proof of Statement 1 of Lemma  $\ref{776}$} 

Let us consider a compact domain in the configuration space, which is presented by the cube with edges of the  length $T$ (the length of the parametrized magnetic line) and decomposes this domain into a finite union of  subdomains.
Subdomains are divided into 3 types as the cases in Lemma  $\ref{77}$. Evidently, in Statement 1,2  the constants $\varepsilon$ could be smaller,
and the corresponding parameter $C$ larger, as required. Finally, we may assume that $\varepsilon$, $\varepsilon_0$, and $C$ are common in the  3 Statements of the lemma.  

Let us construct the corresponding limited tenor, see Definition $\ref{tenzor}$, which estimates the absolute value of the density function
$W_{3,4;2}$ in each subdomain. For a subdomain of the type 1, or 2 we define the limiting tensor as  $C \ln{\delta(x_{1,1},x_{1,2},x_{2,1},x_{3,1})}$.
For a subdomain of the type 3 we define the limiting tensor as $C$.

By the Holder's inequality
we get
$$ \int fg dx \le (\int \left|f\right|^q dx)^{\frac{1}{q}} (\int \left|g\right|^p dx)^{\frac{1}{p}}, \quad 1<p<2, \quad \frac{1}{p} + \frac{1}{q}=1, $$
as in Theorem 5. 
In this inequality $f$ is the limiting tensor for 
$(\ref{jadro})$, which is constructed above using Lemma  $\ref{77}$, $g$--is limiting tensor for 
$(2,3)(3,1)^2(1,2)$, which is constructed analogously as in \cite{A3}, Theorem 5.
We could take arbitrary dominators
 $p,q$, for example,  $p=q=2$. This gives the limiting tensor for $(\ref{jadro})$. Statement 1 is proved.

\subsubsection*{Proof of Statement 2 of Lemma  $\ref{776}$}
Let us start with the case of the 3 closed magnetic lines $L_1,L_2,L_3$ and rewrite the integral $(\ref{37})$
in an "`ergodic style"'. 
Let us construct a configuration space $K$ and calculate the expression
\begin{eqnarray}\label{370}
(3,1)(2,3)D_{1,2} = (3,1)(2,3)\oint \phi^{2}_2  (\dot{x}_2,\A_1) dx_2
\end{eqnarray}
as an integral over this configuration space.

Assume that $(\dot{x}_2,\A_1), \quad x_2 \in L_2$ is a given function and denote this function by $A: L_2 \to \R$.
Take a marked point $pt_2 \in L_2$, take the universal covering $p: \tilde L \to L_2$ over $L_2$ and take a lift $pt \in \tilde L$ of the point $pt_2$.
Denote the periodic function $A \circ p: \tilde L \to \R$ by $\tilde A$.

Denote by $\phi^{pt_2}$ a function $\phi_2: L_2 \to \R$ in the gauge $(\ref{phi2})$, $\phi^{pt_2}(pt_2)=0$. 
Denote by $\tilde{\phi}^{pt}: \tilde L \to \R$ the periodic function $\tilde{\phi}^{pt} = p \circ \phi^{pt_2}$, by the definition 
$\tilde{\phi}^{pt}(pt)=0$.  Denote the period by $l$.

Take the Cartesian product $\tilde L \times \tilde L$ and the diagonal with the marked point $pt \in \tilde L_{diag} \subset \tilde L \times \tilde L$. 
Let us define a coordinate system $(\tau,t)$ in $\tilde L \times \tilde L$, where $\tau$ is the coordinate on the diagonal $\tilde L_{diag}$,
$t$ is the coordinate of the first factor of $\tilde L \times \tilde L$.
Define  the functions $\tilde \phi^{pt}(\tau)$, $\tilde A(\tau)$ as periodic on the diagonal $\tilde L_{diag}$ with the common period 
$l_{diag}=\sqrt{2} \cdot l$, which are induced from the corresponding functions on the first factor by the projection. 
Denote by $\tilde \phi^{pt}_t(\tau)$, the $t$-shift of $\tilde \phi^{pt}$, $\tilde \phi^{pt}_t = \tilde \phi^{pt}(\tau + t)$, $pt, t \in  \tilde L_{diag}$.
Denote by $\psi_{pt,t}: \tilde L_{diag} \to \R$ the  product $\psi_{pt,t} = \tilde \phi^{pt} \cdot \tilde \phi^{pt}_t$. The $t$-family $\psi_{pt,t}$ determines the $(l,l)$-periodic  function $\Psi_{pt}(\tau,t): \tilde L \times \tilde L \to \R$, which is depended of the parameter $pt$.

 For an arbitrary pair $(pt,t)$ of parameters the function $\psi_{pt,t}: \tilde L_{diag} \to \R$ is periodic. Decomposes this function 
as following:  $\psi_{pt,t}(\tau) = m_{pt,t} +  \phi^{var}_{pt,t}(\tau)$, $\tau \in \tilde L$, where  the constant $m_{pt,t}$ is the mean value
of the periodic function $\psi_{pt,t}$. Denote by $m[A]$ the mean value of $A(x)$ and by $m_{t}$ the mean value of $m_{pt,t}$.

Consider the integral of $\tilde \phi^{pt}(t_1 +pt) \cdot \tilde \phi^{pt}(t_2 + pt)$ over the square $K=\{0 \le t_1 \le kl, 0 \le t_2 \le kl\}$, $k \in \N$.
Obviously, we get $\int_{0}^{kl} \tilde \phi^{pt}(pt+t_1)dt_1 = m_{pt,t}kl$, 
$\int \int_K \tilde \phi^{pt}(t_1) \cdot \tilde \phi^{pt}(t_2) dt_1 dt_2=(m_{pt,t}kl)^2$. Obviously, we get the equality:
\begin{eqnarray}\label{KD}
\int \int_K \tilde \phi^{pt}(t_1) \cdot \tilde \phi^{pt}(t_2) dt_1 dt_2 = \int \int_D \Psi_{pt}(\tau,t) d\tau dt, 
\end{eqnarray}
where $D$ is the subdomain in $\tilde L \times \tilde L$, defined by $pt \le \tau \le kl + pt$, $pt \le t \le kl + pt$.

Using this,  rewrite the  function below the integral $(\ref{370})$ in the following asymptotic form:
 \begin{eqnarray}\label{371}
D_{1,2}(x) =  A(x) \lim_{pt \to -\infty} \frac{d(pt,x)}{(pt-x)^2}, \quad pt, x \in \tilde L_2,
\end{eqnarray}
$$  d(pt,x) = \int_{pt}^{x} \int_{pt}^{x} \Psi_{pt}(t,\tau)  d t d \tau.$$

To simplify this equation take the mean value $\psi_{pt,t}(\tau) \mapsto m_{t}$, where $pt \to -\infty$. Take the mean value $A(x) \mapsto (1,2)$ over $x \in \tilde L_2$, we get:
  \begin{eqnarray}\label{372}
D_{1,2} = (1,2) \lim_{t_0 \to +\infty} \frac{\int_0^{t_0}  m(t) dt}{t_0}.  
\end{eqnarray}


For an arbitrary $t \in [0,+\infty)$ define the configuration space $K^t=K^t_{1,2,1}$ and calculate the value of the function $m(t)$  in the formula $(\ref{372})$ for a given $t$ as the asymptotic integral over $K^t$. 
The integral $(\ref{37})$ is generalized for generic magnetic fields $\B \in \Omega$ and the function $m(t)$ is bounded and determines the asymptotic functional (generally speaking, this functional is multivalued). 

Assume  that the linking numbers $(1,2),(2,3),(3,1)$ are fixed.
Define the function $m(t)$ by Theorem $\ref{th4}$ for $j=1$ as the asymptotic integral over $1+2+1$-points configuration space $K^t_{1,2,1}$, where $t$ is the
the time of the magnetic flow, which transforms the points $y_1 \in \tilde L_2$ into the point $y_2 \in \tilde L_2$. 

To calculate the term $(\ref{371})$ we apply Theorem $\ref{th4}$ for $j=4$, where the factors of classifying space for  $j=2,3,4$  correspond
with the asymptotic linking numbers $(1,2),(2,3),(3,1)$.

The value $m(t)$ for an arbitrary $t \ge 0$ is well-defined and the absolute value $\vert m(t) \vert$ is bounded  for an arbitrary $t$ by a constant, which depends no of $t$. The required integral $(\ref{371})$ is defined by the upper and lower asymptotic limits in the formula $(\ref{372})$.

Lemma $\ref{776}$ is proved.

\[  \]

The following statement is the main theorem of the paper.

\begin{theorem}\label{mainTh}

1. There exists the asymptotic functional  $M$ on the space
$\Omega$ in the sense of Definition $\ref{local}$, given by the sum as in the formula $(\ref{40})$. The main term   $(\ref{18})$
is well defined, and the term $(\ref{37})$ is well-defined as a multivalued real between lower and upper bounds.

2. Assume that the magnetic field  $\B \in \Omega$ satisfies Definition $\ref{vosvr}$. 
Then the asymptotic functional  $M$, which is construction in 1, coincides with the integral of values of the combinatorial invariant 
$(\ref{40})$ over arbitrary ordered triples of (closed) magnetic lines of  $\B$ (the measure on the space of closed magnetic lines is defined
by  uniformity of the integral magnetic flow over cross-sections of magnetic lines). The local functional $M$ is an invariant of volume-preserved transformation of the space.

3. Assume that the magnetic field $\B \in \Omega$ is decomposed into 3 disjoint magnetic tubes. Then the asymptotic functional  $M$, which is constructed in 1., which is  defined by the integration of expressions 
$(\ref{40})$ over arbitrary ordered triples of (generally speaking, non-closed) magnetic lines of  $\B$, each magnetic line is contained inside the corresponding magnetic tube, coincides with the origin integral invariant $(\ref{MM})$.

\end{theorem}

\subsubsection*{Proof of Theorem  $\ref{mainTh}$}

Statement 1 follows from Lemma 
$\ref{776}$. 

Statement 2 follows from Theorem 
$\ref{th1}$, Statement 2. (Let us remind that in the formula $(\ref{40})$ the lengths
of parameterizations of magnetic lines are arbitrary. 

Statement 3 is obvious by the construction.
Theorem $\ref{mainTh}$ is proved. 
\[  \]

\section{Conclusion Remarks}
The goal of the paper is to introduce a new invariant $M$ for magnetic fields, which is written in an "`ergodic style"', analogously to the asymptotic ergodic Hopf invariant. Statement 1 of Main Theorem is proved only for terms $(\ref{18})$, $(\ref{37})$. The authors clams that all the terms in the integral $(\ref{40})$ the proof is analogous. The question if some of the terms are multivalued is open.
The question if the integral, introduced in Theorem $\ref{mainTh}$, determines an invariant for generic ideal magnetic fields up to volume-preserved diffeomorphisms is open.
The author conjecture the answer is affirmative.

\end{document}